\documentclass[aos,MSNbibl,seceqn,citesort,dvips]{arximspdf}
\usepackage{mathrsfs}
\usepackage{graphicx}
\usepackage{vtexurl}

\doi{10.1214/11-AOS954}
\volume{40}
\issue{1}
\pubyear{2012}
\firstpage{159}
\lastpage{187}

\makeatletter

\newcommand{\iint}{\int\!\!\int}

\newtheorem{theorem}{Theorem}
\newtheorem{lemma}{Lemma}

\newproclaim{remark}{Remark}
\newtheorem{claim}{Claim}

\newcommand{\argmin}{\mathop{\arg\min}}

\def\bsuffix #1{#1}

\setattribute{abstract}{width}{280pt}
\makeatother

\begin{document}
\begin{frontmatter}

\title{Large-sample study of the kernel density estimators under
multiplicative censoring\thanksref{T3}}
\runtitle{KDE under multiplicative censoring}

\begin{aug}
\author[A]{\fnms{Masoud} \snm{Asgharian}\corref{}\thanksref{t1}\ead[label=e1]{masoud@math.mcgill.ca}},
\author[B]{\fnms{Marco} \snm{Carone}\thanksref{t2}\ead[label=e2]{mcarone@berkeley.edu}}
\and
\author[C]{\fnms{Vahid} \snm{Fakoor}\ead[label=e3]{fakoor@math.um.ac.ir}}
\runauthor{M. Asgharian, M. Carone and V. Fakoor}
\affiliation{McGill University, University of California,
Berkeley, and~Ferdowsi~University~of~Mashhad}
\address[A]{M. Asgharian\\
Department of Mathematics and Statistics\\
McGill University\\
Burnside Hall, Room 1237\\
805 rue Sherbrooke Ouest\\
Montr\'{e}al, QC, H3A 2K6\\
Canada\\
\printead{e1}}
\address[B]{M. Carone\\
Division of Biostatistics\\
University of California, Berkeley\\
109 Haviland Hall\\
Berkeley, California 94720-7358\\
USA \\
\printead{e2}}
\address[C]{V. Fakoor\\
Department of Statistics\\
Ferdowsi University of Mashhad\\
Azadi Square\\
Mashhad, Razavi Khorasan 1159-91775\\
Iran\\
\printead{e3}} 
\end{aug}

\thankstext{T3}{Authorship is listed alphabetically to reflect the equal
contribution of each author.}

\thankstext{t1}{Supported in part by NSERC, FQRNT and the
Veterans Affairs HSR\&D Service.}

\thankstext{t2}{Supported in part by NSERC, FQRNT, NIH and the Hopkins Sommer
Scholars Program.}

\received{\smonth{5} \syear{2009}}
\revised{\smonth{7} \syear{2011}}

%
\begin{abstract}
The multiplicative censoring model introduced in Vardi
[\textit{Biomet\-rika} \textbf{76} (1989) 751--761] is an incomplete data
problem whereby two independent samples from the lifetime distribution
$G$, $\mathcal{X}_m=(X_1,\ldots,X_m)$ and
$\mathcal{Z}_n=(Z_1,\ldots,Z_n)$, are observed subject to a form of
coarsening. Specifically, sample $\mathcal{X}_m$ is fully observed
while $\mathcal{Y}_n=(Y_1,\ldots,Y_n)$ is observed instead of
$\mathcal{Z}_n$, where $Y_i=U_iZ_i$ and $(U_1,\ldots,U_n)$ is an
independent sample from the standard uniform distribution. Vardi
[\textit{Biometrika} \textbf{76} (1989) 751--761] showed that this
model unifies several important statistical problems, such as the
deconvolution of an exponential random variable, estimation under a
decreasing density constraint and an estimation problem in renewal
processes. In this paper, we establish the large-sample properties of
kernel density estimators under the multiplicative censoring model. We
first construct a~strong approximation for the process
$\sqrt{k}(\hat{G}-G)$, where $\hat{G}$ is a solution of the
nonparametric score equation based on $(\mathcal{X}_m,\mathcal{Y}_n)$,
and $k=m+n$ is the total sample size. Using this strong approximation
and a result on the global modulus of continuity, we establish
conditions for the strong uniform consistency of kernel density
estimators. We also make use of this strong approximation to study the
weak convergence and integrated squared error properties of these
estimators. We conclude by extending our results to the setting of
length-biased sampling.
\end{abstract}

%

\end{frontmatter}

\section{Introduction} \label{intro}

Vardi~\cite{vardi1989biometrika} introduced an incomplete data problem
unifying several statistical models. The problem consisted of inferring
the lifetime distribution of interest $G$ through a random sample
$X_{1}, X_2,\ldots, X_{m}$ drawn directly from $G$ and a random sample
$Y_1,Y_2,\ldots,Y_n$ drawn from the distribution $F$ with density function
%
\begin{equation} \label{resid-life}
f(y)=\int_{y \leq z} z^{-1}\,d G(z),\qquad y > 0.
\end{equation}
Since $f$ is a decreasing density function, $Y$ may be expressed as the
product of two independent random variables: a nonnegative variate $Z$
and a standard uniform variate $U$. From the form of~(\ref
{resid-life}), it is easy to see that in this case $Z$ must be
distributed according to $G$. This representation suggests that only a
random fraction of $Z$ may be observed, motivating the nomenclature
\textit{multiplicative censoring} used to describe this incomplete data
scheme. The likelihood based on the $k=m+n$ observations $X_{1}=x_1,
\ldots, X_{m}=x_m$ and $Y_{1}=y_1, \ldots, Y_{n}=y_n$ is
%
\begin{equation} \label{mc-like}
L(G)=\prod_{i=1}^{m} G(dx_{i}) \prod_{j=1}^{n} \int_{y_{j} \leq z}
z^{-1}\,dG(z).
\end{equation}
As discussed by Vardi~\cite{vardi1989biometrika}, the multiplicative
censoring model arises from the deconvolution of an exponential random
variable, estimation under a decreasing density constraint and an
estimation problem in renewal processes. The literature on these and
related problems is vast. Estimation under a~decreasing density
constraint dates back to the seminal work of Grenander~\cite
{grenander1956actscand}, with key contributions by Groeneboom~\cite
{groeneboom1985berkconf} and Huang and Wellner~\cite{huang1995scand}.
The estimation problem in renewal processes discussed in~\cite
{vardi1989biometrika} is closely tied to important applications in
cross-sectional sampling and prevalent cohort studies in epidemiology
(length-biased sampling) and in labor force studies in economics
(stock sampling). The multiplicative censoring model and its variants
have been studied by
\cite{hasminskii1983ln,vardi1989biometrika,bickel1994theoryprob,bickel1993,vandervaart1994annals}
and~\cite{vardi1992lannals}, among others. Vardi~\cite{vardi1992lannals}
studied the asymptotic behavior of solutions of the nonparametric score
equation under the multiplicative censoring model.

As will be discussed later, multiplicative censoring and left-truncated
right-censored data are intricately tied. The latter have been
extensively studied in the statistical literature. Their importance
stems mainly, although not exclusively, from the widespread use of
prevalent cohort study designs to estimate survival from onset of a
disease. In such studies, patients with prevalent disease are
identified at some instant in calendar time, often through a
cross-sectional survey. These patients are then followed forward in
time until death or loss to follow-up. If no temporal change in the
incidence of disease has occurred during the period covering observed
onsets, a stationary Poisson process may adequately describe the
incidence pattern of the disease; see
\cite{asgharian2006statsinmed,asgharian2002jasa,asgharian2005annals} and
\cite{wolfson2001nejm}. In this case, the left-truncation variable
is uniformly distributed, and the failure time data are said to be
length-biased. The likelihood for the observed data is then given by
(\ref{mc-like}), where
\[
G(t)=\mu_U^{-1}\int_{0}^{t}{u}\,dF_{U}(u),
\]
$\mu_U=\int_{0}^{\infty}{u}\,dF_{U}(u)$ and $F_{U}$, the \textit{unbiased
distribution}, is the underlying distribution function about which we
would like to infer; see Section~\ref{lbs-rc} and~\cite
{asgharian2005annals}. Because we require $\mu_U<\infty$ in the above,
we restrict our attention to distribution functions $G$ such that $\int
{z^{-1}\,dG(z)}<\infty$.

%

The connection between the multiplicative censoring model and prevalent
cohort studies under the stationarity assumption has revived interest
in the former. Nonetheless, there appears to be no result in the
literature on density estimation under the multiplicative censoring
model, despite its importance in applied sciences. A recent application
described by Kvam~\cite{kvam2008tech} concerns nanoscience and the
measurement of carbon nanotubes. As discussed by Silverman~\cite
{silverman1986}, density estimation can be useful for purposes of data
exploration and presentation. It is effective in the investigation of
modes (determination of multimodality and identification of modes) and
tail behavior (rate of tail decay). These features are especially
important in length-biased sampling and survival analysis, where
skewness is often pervasive and differential subgroup characteristics
may lead to multimodality. An additional motivation for the study of
density estimation under multiplicative censoring stems from the fact
that nonparametric regression of right-censored length-biased data has
not been addressed in the literature. In view of the intricate link
between density estimation and nonparametric regression (see~\cite
{nussbaum1996annals}), a~study of density estimation under
multiplicative censoring provides foundations for studying
nonparametric regression of right-censored length-biased data.

Among the various methods of density estimation, kernel smoothing is
particularly appealing for both its simplicity and its interpretability
(e.g., as a limiting pointwise average of shifted histograms). It
provides a unifying framework in that, as discussed in~\cite
{scott1992}, each of finite difference density estimation, smoothing by
convolution, orthogonal series approximations and other smoothing
methods historically used in the various applied sciences can be seen
as instances of kernel smoothing. This article studies the large-sample
properties of kernel density estimators in the setting of
multiplicative censoring. Pioneered by Silverman~\cite
{silverman1978annals}, the approach adopted consists of constructing
strong approximations of the empirical density process.

Although under the multiplicative censoring model we may avoid
complexities altogether by performing estimation using the uncensored
observations alone, use of the full data is motivated by at least two
reasons. First, although discarding the censored cases under the
canonical multiplicative censoring scheme does not compromise
consistency, the same cannot be said under the related length-biased
sampling scheme, even though these schemes lead to the same likelihood.
This occurs because, under length-bias sampling, the uncensored cases
do not emanate directly from the (length-biased version of the)
distribution of interest. Systematic exclusion of the censored cases
would therefore lead to inconsistency. This fact motivates the study of
both censored and uncensored cases under multiplicative censoring.
Second, due to the informativeness of the censoring mechanism, ignoring
the censored observations may lead to a substantial loss of efficiency.
Because the asymptotic covariance
function of the nonparametric maximum likelihood estimator of $G$ does
not have an explicit form, this phenomenon is difficult to quantify in
the nonparametric setting (see the discussion on page 1024 of~\cite
{vardi1992lannals}); however, a~parametric example may be illustrative.
Suppose that the uncensored observations emanate from a Gamma
distribution, say with mean 2$\theta$ and variance~$2\theta^2$, then
the censored observations are exponentially distributed with mean
$\theta$. The asymptotic relative efficiency of the full-sample MLE
relative to the uncensored-sample MLE is $1+\upsilon/2$, where
$\upsilon
>0$ is the asymptotic relative frequency of censored observations to
uncensored observations. If, for example, $\upsilon=1$, indicating that
uncensored and censored cases arise in equal numbers asymptotically,
use of the full sample provides a fifty percent gain in efficiency.

Following~\cite{komlos1975zwvg}, hereafter referred to as KMT, and
\cite
{csorgo1981}, we first construct a strong approximation for the process
$\sqrt{k} (\hat{G} - G)$, where $\hat{G}$ is a solution of the
nonparametric score equation based on $(\mathcal{X}_m,\mathcal{Y}_n)$.
The literature on strong approximations is vast. Recent reviews on
empirical processes, strong approximations and the KMT construction
include~\cite{delbarrio2007} and~\cite{mason2007jspi}. Using this
strong approximation and a result on the global modulus of continuity,
we obtain the strong uniform consistency of the kernel density
estimators of the density function $g$ associated to $G$ and find a
sequence of Gaussian processes strongly uniformly approximating the
empirical kernel density process. Using these results, we study the
integrated squared error properties of the kernel density estimators.

The layout of the paper is as follows. In Section~\ref{prelim}, we
introduce our notation and present some preliminaries. In Section~\ref
{st-ap}, we find a sequence of Gaussian processes that strongly
uniformly approximates the empirical process $\sqrt{k} (\hat{G} - G)$
and study its global modulus of continuity. We use these results to
study the asymptotic behavior of the kernel density estimators in
Section~\ref{kernel}. It is shown, in particular, that the kernel
density estimators are strongly consistent and asymptotically Gaussian.
Section~\ref{ISE} is devoted to the integrated squared error properties
of the kernel density estimators and includes results from a
preliminary small-sample simulation study. We show how our results can
be extended to length-biased sampling with right-censoring in Section
\ref{lbs-rc} and present concluding remarks in Section~\ref
{conclusion}. The claim and theorems are proved in the
\hyperref[app]{Appendix}
while lemmas are proved in the supplementary material~\cite
{asgharian2012supp}.


\section{Preliminaries} \label{prelim}
We consider the random multiplicative censoring model introduced in
\cite{vardi1989biometrika}, whereby two independent random samples
$\mathcal{X}_m=(X_1,\ldots,\allowbreak X_m)$ and $\mathcal{Z}_n=(Z_1,\ldots,Z_n)$
are drawn from the lifetime distribution $G$ and a~third independent
sample $\mathcal{U}_n=(U_1,\ldots,U_n)$, from the standard uniform
distribution. Let $Y_i=Z_iU_i$, $i=1,\ldots,n$, and write $\mathcal
{Y}_n=(Y_1,\ldots,Y_n)$. Then~$\mathcal{Y}_n$ is a random sample from
the absolutely continuous distribution~$F$ with density given by~(\ref
{resid-life}). The observed data consist of $(\mathcal{X}_m,\mathcal
{Y}_n)$ while $(\mathcal{Z}_n,\mathcal{U}_n)$
is unobserved.

We begin with the score equation derived from the likelihood $L(G)$
given by~(\ref{mc-like}). Let $G_m$ and $F_n$ be, respectively, the
empirical distribution functions based on the uncensored observations
$x_1,\ldots,x_m$ and the censored cases $y_1,\ldots,y_n$, and write
$\hat{p} = m/k$, where $k = m+n$. For simplicity, assume all
observations are distinct, and denote by $t_{1} < \cdots<t_{k}$ the
values taken by $x_{1}, \ldots, x_{m}$ and $y_{1}, \ldots,
y_{n}$. The distribution function $\hat{G}$ satisfies the nonparametric
score equation if, for all $t\geq0$,
%
\begin{equation} \label{sq1}
d\hat{G}(t)=\hat{p} \,dG_{m}(t) + (1-\hat{p}) \biggl[\int_{0 < y \leq t}
\frac{dF_{n}(y)}
{\int_{y \leq z} z^{-1}\,d\hat{G}(z)} \biggr] t^{-1}
\,d\hat{G}(t),\hspace*{-25pt}
\end{equation}
while $\sum_{j=1}^{k} d\hat{G}(t_{j}) = 1$ and $d\hat{G}(t_{j}) \geq0$,
$j=1, \ldots, k$; see~\cite{vardi1992lannals}, page 1025. Integrating both
sides of~(\ref{sq1}), we obtain
\[
\hat{G}(t)=\hat{p} G_{m}(t) + (1-\hat{p}) \int_{0 <x \leq t}
\biggl[\int_{0 < y \leq x} \frac{dF_{n}(y)} {\int_{y \leq z}z^{-1} d\hat{G}(z)}
\biggr] x^{-1} \,d\hat{G}(x),
\]
where the final integrand is defined to be 0 for $x > t_{k}$. We say
that a sequence of real numbers $\gamma_{m,n}$ satisfies assumption
(A0) if
\[
\sum_{m,n}G(\gamma_{m,n})<\infty,
\]
where the summation is understood to range over subsample sizes $m$ and~$n$,
jointly taken to infinity, so that $\hat{p}\rightarrow p\in(0,1]$.
To circumvent problems related to a singularity at the origin, we
select a sequence of positive real numbers $\gamma_{m,n}$ satisfying
(A0) and consider solutions $\hat{G}$ of~(\ref{sq1}) assigning zero
mass below $\gamma_{m,n}$. All results derived in this article apply to
any solution of~(\ref{sq1}) with this property. The existence of such
solutions is an important fact.
%
%
\begin{claim}\label{cla1}
Suppose that \textup{(A0)} holds. Then, for each $m$ and $n$ sufficiently large,
(\ref{sq1}) has a solution $\hat{G}$ such that $\hat{G}(u)=0$ for each
$u<\gamma_{m,n}$.
\end{claim}

If there exists some $\gamma_0>0$ such that $G(\gamma_0)=0$, assumption
(A0) is not required. We may simply choose $\gamma_{m,n}=\gamma_0$, and
because any solution of~(\ref{sq1}) will have zero mass below $\gamma
_{m,n}$, the proposition follows directly from~\cite{vardi1989biometrika}.


Define
$U_{m,n} = \sqrt{k}(\hat{G}-G)$, $W_{X,m} = \sqrt{m}(G_m-G)$,
$W_{Y,n} = \sqrt{n}(F_n-F)$, $\hat{f}(t)=\int_{t\leq z}z^{-1}\,d\hat
{G}(z)$ and
%
\begin{equation}\label{e4}\qquad
W_{m,n}(t) = \sqrt{\hat{p}} W_{X,m}(t)+\sqrt{1-\hat{p}} \hat{f}(t)
\int_{0<y\leq t}W_{Y,n}(y)\,d\biggl[\frac{1} {\hat{f}(y)} \biggr].
\end{equation}
We observe, in particular, that
%
\begin{equation}\label{boundwmn}
|W_{m,n}(t)|\leq\sqrt{\hat{p}}|W_{X,m}(t)
|+\sqrt
{1-\hat{p}}\sup_{0<y\leq t}|W_{Y,n}(y)|
\end{equation}
for each $t>0$.
As in~\cite{vardi1992lannals}, we have that
\[
W_{m,n}(t) = \hat{p} U_{m,n}(t)+(1-\hat{p})\hat{f}(t)\int_{0<y \leq
t } y
\biggl(\int_{y \leq z}\frac{U_{m,n}(z)}{z^{2}} \,dz \biggr)\,
d\biggl[\frac{1}{\hat{f}(y)} \biggr].
\]
The process $W_{m,n}$ can therefore be expressed as the image of a linear
operator applied on $U_{m,n}$.
To see this, we define the operator $\mathcal{G}_{m,n}$ pointwise as
$\mathcal{G}_{m,n}(u)(t)
= \hat{f}(t)\mathcal{A}_{m,n}(u)(t)$, where
\[
\mathcal{A}_{m,n}(u)(t)=\int_{0<y \leq t }y
\biggl(\int_{y \leq z}\frac{u(z)}{z^{2}} \,dz \biggr) \,d\biggl[
\frac{1}{\hat{f}(y)} \biggr].
\]
Then, we may write $\mathcal{F}_{m,n}=\hat{p}\mathcal{I}+(1-\hat
{p})\mathcal{G}_{m,n}$, with $
\mathcal{I}(u) = u$ the identity map,
and observe that
%
\begin{equation} \label{linear}
W_{m,n} = \mathcal{F}_{m,n}(U_{m,n}).
\end{equation}
Denoting by $D_{0} [0, \infty]$ the space of \textit{cadlag} functions
vanishing at 0 and $\infty$ endowed with the uniform
topology (the topology induced by the supremum norm over $[0,\infty)$,
$\| u
\|_{\infty} = {\sup_{0 \leq t < \infty}} |u(t)|$), it is not difficult to
see that $\mathcal{I}$, $\mathcal{G}_{m,n}$ and $\mathcal{F}_{m,n}$ are
bounded linear operators on $D_{0} [0, \infty]$, and, in view of Lemma
3 of~\cite{vardi1992lannals}, that $\mathcal{F}_{m,n}$ has a bounded
inverse satisfying $\|\mathcal{F}_{m,n}^{-1}\| \leq2/\hat{p}^{2}$. As
in~\cite{vardi1992lannals}, it holds that if $\hat{p}\rightarrow p\in
(0,1]$ as $m, n \rightarrow\infty$, then, for each $u\in D_0[0,\infty
]$, we have that
\[
\| \mathcal{F}_{m,n}(u) - \mathcal{F}(u) \|_{\infty}
\stackrel{\mathrm{a.s.}}{\longrightarrow} 0,
\]
where the limit operators are $\mathcal{F} = p \mathcal{I} +
(1-p)\mathcal{G}$, $\mathcal{G}(u)(t)=f(t)\mathcal{A}(u)(t)$ and
\[
\mathcal{A}(u)(t)=\int_{0<y \leq t }y \biggl(\int_{y
\leq z} \frac{u(z)}{z^{2}} \,dz \biggr) \,d
\biggl[\frac{1}{f(y)} \biggr].
\]
We may then conclude that $\mathcal{G}$ and $\mathcal{F}$ are also
bounded linear operators on $D_{0} [0, \infty]$ and that $\mathcal{F}$
has a bounded inverse satisfying $\| \mathcal{F}^{-1}\|
\leq2/p^{2}$. Vardi~\cite{vardi1992lannals}
proved the uniform strong consistency of $\hat{G}$ using~(\ref
{linear}). Instead, we obtain it as a corollary of Lemma~\ref{lemma31} below.

Of importance will be the fact, proved in~\cite{vardi1992lannals}, that
the inverse operator~$\mathcal{F}^{-1}$ has the following pointwise
representation:
%
\begin{equation} \label{represent}
\mathcal{F}^{-1}(u)(t)=p^{-1}u(t)+\int_{0}^{\infty}{\Phi
(t,x)u(x)\,dx}
\end{equation}
with kernel $\Phi$ satisfying, for each $t$ and $x$, the constraints
%
\begin{equation} \label{constraint1}\quad
p^2\Phi(t,x)+(1-p)\mathcal
{A}_0(t,x)+p(1-p)\int_{0}^{\infty}{\Phi(t,z)\mathcal{A}_0(z,x)\,dz} = 0
\end{equation}
and
%
\begin{equation} \label{constraint2}\int_{0}^{\infty
}{\Phi(t,z)\mathcal{A}_0(z,x)\,dz}= \int_{0}^{\infty}{\mathcal
{A}_0(t,z)\Phi(z,x)\,dz},
\end{equation}
where we have defined $\mathcal{A}_0(t,x)=f(t)x^{-2}\int_{0<y\leq
t\wedge x}{y\,d[1/f(y)]}$.\vspace*{2pt}

As in~\cite{vardi1992lannals}, we have that $W_{m, n}
\leadsto W$ in $D_{0}[0, \infty]$, where $W$ is the
Gaussian process
\[
W(t) = \sqrt{p} B_{X} (G(t)) + \sqrt{1-p} f(t) \int_{0 < y \leq t}
B_{Y} (F(y)) \,d\biggl[\frac{1}{f(y)} \biggr]
\]
with $B_X$ and $B_Y$ independent Brownian bridges, and that $U_{m,n}
\leadsto U= \mathcal{F}^{-1} (W)$ in $D_{0}
[0, \infty]$. Here, the symbol $\leadsto$ refers to weak convergence.
This last step can be established using the convergence of
$\mathcal{F}_{m,n}$ to $\mathcal{F}$ in operator norm topology,
Lemma 3 of~\cite{vardi1992lannals} and the continuous mapping
theorem. A~consistent estimator $\hat{\psi}_U(s,t)$ of $\psi
_U(s,t)=\mathrm{E}[U(s)U(t)]$ is provided in~\cite
{vardi1992lannals}, though in practice the use of resampling methods
may yield an estimator of~$\psi(s,t)$ more expediently.


\section{Approximation of the empirical process $U_{m,n}$}\label{st-ap}

\subsection{Strong approximation}

Let $\alpha_n$ denote the empirical process of $n$ independent
standard uniform random variables. The KMT construction implies that
there exists a
probability space $(\Omega,\mathcal{F},P)$ with a sequence of
independent standard uniform random variables and a sequence of Brownian
bridges $B_n$ such that
\[
\|\alpha_n-B_n\|_{[0,1]} = \mathcal{O}\biggl(\frac{\log n}{\sqrt
{n}}\biggr) \qquad\mbox{a.s.}
\]
Equation~(\ref{linear}) is key to the strong approximation of
$U_{m,n}$. Since $W_{X,m}$ and~$W_{Y,n}$ are
independent empirical processes associated, respectively,\break with~$\mathcal
{X}_m$ and $\mathcal{Y}_n$, in view of the KMT construction, there
exist versions of~$W_{X,m}$ and $W_{Y,n}$ along with two independent sequences of Brownian
bridge processes~$B_{X,m}$ and
$B_{Y,n}$ such that $B_{X,m}\circ G$ and $B_{Y,n}\circ F$ approximate~$W_{X,m}$
and $W_{Y,n}$ at the optimal rate of $\log s / \sqrt{s}$ (here,
$s$ is the sample size). Using~(\ref{linear}), we extend this
approximation to $W_{m,n}$ and use properties of $\mathcal{F}$
to find a sequence of Gaussian processes strongly uniformly
approximating~$U_{m,n}$. The main theorem of
this section, Theorem~\ref{theorem31}, is proved through a~sequence of
lemmas.

Denote the upper limit of the support of $G$ by $\tau=\sup\{t\dvtx G(t)<1\}
$. Given any set $B$, denote by $\mathbb{I}_B$ and \mbox{$\|\cdot\|_B$} the
indicator function and the supremum norm over~$B$, respectively. Write
\mbox{$\|\cdot\|_\infty$} for the case $B=[0,\infty)$. We introduce the
following assumptions:\vspace*{8pt}

(A1) $\sqrt{k}(\hat{p} - p)
= \mathcal{O}(\sqrt{\log\log k} )$ for some $p \in
(0,1]$.

(A2) $G$ is continuous and has bounded support ($\tau<\infty
$).

(A3) There exists $\alpha_0>2$ such that $\lim_{x\downarrow
0}G(x)/x^{\alpha_0}<\infty$.

(A4) There exists $\beta>0$ such that $\lim_{x\downarrow
0}
[1-G(\tau-x)]/x^{\beta}\in(0,\infty)$.\vspace*{8pt}

We begin\vspace*{1pt} by obtaining rates for the difference between $\hat{G}$ and
$G$ as well as between $\hat{f}$ and $f$ in the supremum norm.
%
\begin{lemma}\label{lemma31} Suppose \textup{(A0)} holds. Then, for any
sequence of nonnegative real numbers $a_{m,n}$, as $k\rightarrow
\infty$:
\begin{eqnarray*}
&&\mbox{\textup{(a)}}\quad \| \hat{G}-G \|_{\infty} = \mathcal{O}\Biggl(\sqrt{\frac{\log
\log
k}{k}}\Biggr) \qquad\mbox{a.s.},\\
&&\mbox{\textup{(b)}}\quad \| \hat{f}-f \|_{[a_{m,n},\infty)} \\
&&\hphantom{\mbox{\textup{(b)}}}\qquad
=\mathcal{O}\Biggl(\gamma_{m,n}^{-1}\sqrt{\frac{\log
\log
k}{k}}+[F_U(\gamma_{m,n})-F_U(a_{m,n})]\mathbb
{I}_{[0,\gamma
_{m,n})}(a_{m,n})\Biggr) \quad\mbox{a.s.}
\end{eqnarray*}
\end{lemma}

The above indicates, for example, that in addition to satisfying (A0),
$\gamma_{m,n}$ should be such that
\[
\gamma_{m,n}^{-1}\sqrt{\frac{\log\log k}{k}}\rightarrow0.
\]
If (A3) holds, the sequence $\gamma_{m,n}'=k^{-{1}/({2\alpha})}$ may
be considered, with the choice $\alpha\in(1,\alpha_0/2)$ ensuring that
the two requirements above are satisfied. In this case, choosing
$\alpha
$ as close as possible to $\alpha_0/2$ would yield the fastest rate,
modulo logarithmic terms, in part (b) of Lemma~\ref{lemma31}. We now provide
a~result on the growth rate of maxima of Wiener processes.
%
\begin{lemma}\label{lemma32} Let $\mathcal{W}_n$
be a sequence of standard Wiener processes. Then, as $n\rightarrow
\infty$,
\[
\|\mathcal{W}_n\|_{[0,1]}=\mathcal{O}\bigl(\sqrt{\log n}\bigr)
\qquad\mbox{a.s.}
\]
\end{lemma}

The next result considers the asymptotic behavior of the sequence of
inverse operators $\mathcal{F}^{-1}_{m,n}$. First, we note that the
space $D_{0}[0, \tau]$ endowed with the uniform topology is a Banach
space. As such, $\mathscr{A} =
\mathcal{L}(D_{0}[0,\tau], D_{0}[0,\tau])$, the space of bounded
linear operators on $D_{0}[0,\tau]$ endowed with the operator norm topology,
is a Banach algebra. We recall additionally that \textit{cadlag}
functions have countably many jumps (see~\cite{parthasarathy2005}) and
are therefore Riemann integrable on bounded intervals.

Fixing $\varepsilon>0$, set $\mathcal{I}_{\varepsilon}(u)(t)=u(t)\mathbb
{I}_{[0,\tau-\varepsilon]}(t)$ and define $\mathcal{F}_{m,n, \varepsilon}$
and $\mathcal{F}_{\varepsilon}\dvtx D_{0}[0$, $\tau] \to D_{0}[0, \tau]$ as
\[
\mathcal{F}_{m,n, \varepsilon} = \hat{p} \mathcal{I} + (1-\hat{p})
\mathcal
{G}_{m,n,\varepsilon}
\quad\mbox{and}\quad
\mathcal{F}_{\varepsilon} = p \mathcal{I} + (1-p) \mathcal
{G}_{\varepsilon},
\]
respectively, where for any $t \in[0, \tau]$,
\[
\mathcal{G}_{m,n, \varepsilon}(u)(t) = \hat{f}(t) (\mathcal
{A}_{m,n}\circ\mathcal{I}_{\varepsilon})(u)(t)
\quad\mbox{and}\quad
\mathcal{G}_{\varepsilon}(u)(t)
= f(t) ( \mathcal{A}\circ\mathcal{I}_{\varepsilon} )
(u)(t).
\]
Define $\varepsilon_0=\tau p^2/(p^2-2p+2)$.
%
\begin{lemma}\label{lemma39} Suppose that
\textup{(A0)--(A2)} hold and that
$\varepsilon$ is in $(0, \varepsilon_0)$. Then, considering the operator norm
over the space $C_0[0,\tau]$ of continuous functions on $[0,\tau]$
vanishing at the endpoints, as $k\rightarrow\infty$,
\begin{eqnarray*}
&&
\| \mathcal{F}_{m,n, \varepsilon}^{-1} - \mathcal{F}_{\varepsilon}^{-1}\|
\\
&&\qquad=
\mathcal{O}\Biggl(\biggl[\frac{\log(1/\gamma_{m,n})}{f(\tau-
\varepsilon)}
+\frac{F_{U}(\gamma_{m,n})}{f(\gamma_{m,n})}\biggr]\gamma
_{m,n}^{-1}\sqrt{\frac{\log\log k}{k}} + F_U(\gamma_{m,n})
\Biggr) \qquad\mbox{a.s.}
\end{eqnarray*}
%
\end{lemma}

With the choice $\gamma_{m,n}=\gamma_{m,n}'$, the order above may be
simplified to
\[
\|\mathcal{F}_{m,n, \varepsilon}^{-1} - \mathcal{F}_{\varepsilon}^{-1}
\| =
\mathcal{O}\biggl(\frac{k^{-({\alpha-1})/({2\alpha})}\log k\sqrt
{\log
\log k}}{f(\tau-\varepsilon)}\biggr) \qquad\mbox{a.s.}
\]
We now consider a random integral useful in determining the rate of the
strong approximation we will construct for $U_{m,n}$.
%
\begin{lemma}\label{new4} Suppose that \textup{(A0)--(A2)} hold and that
$\varepsilon$ is in $(0, \varepsilon_0)$. Then, as $k\rightarrow\infty$,
\begin{eqnarray*}
&&
\sup_{0\leq s \leq\tau- \varepsilon} \hat{f}(s) \biggl| \int
_{0}^{s}{B_{Y,n}(F(y))\,d\biggl[\frac{1}{\hat{f}(y)} -\frac
{1}{f(y)}
\biggr]}\biggr|
\\
&&\qquad=\mathcal{O}\biggl(\frac{k^{-{1}/{4}}\sqrt{\log
k}(\log
\log k)^{{1}/{4}}}{f(\tau-\varepsilon)}\biggr) \qquad\mbox{a.s.}
\end{eqnarray*}
\end{lemma}
%
\begin{remark}
The above bound also holds for $\varepsilon=\varepsilon_{m,n}\downarrow0$
provided $\varepsilon_{m,n}k/\allowbreak\sqrt{\log\log k}\rightarrow\infty$.
\end{remark}

Henceforth, we set $\gamma_{m,n}=\gamma_{m,n}'$ for each $m$ and $n$.
The next lemma establishes the existence of a sequence of Gaussian
processes approximating~$W_{m,n}$. Define the sequence of processes
%
\begin{equation}\label{32}
W_{m,n}^0(s) = \sqrt{p}B_{X,m}(G(s))+\sqrt{1-p}f(s)\int_{0< y \leq
s}B_{Y,n}(F(y))\,d\biggl[\frac{1}{f(y)}\biggr] .\hspace*{-35pt}
\end{equation}

\begin{lemma}\label{lemma35} Suppose that \textup{(A1)--(A3)} hold and that
$\varepsilon$ is in $(0, \varepsilon_0)$. Then, setting there exists a
probability space on which $W_{m,n}$ and $W^0_{m,n}$ are defined such
that, as $k \to\infty$,
\[
\|W_{m,n}-W_{m,n}^0\|_{[0,\tau-\varepsilon]}=
\mathcal{O}\biggl(\frac{k^{-r(\alpha)}\sqrt{\log k}(\log\log
k)^{{1}/{4}}}{f(\tau-\varepsilon)}\biggr) \qquad\mbox{a.s.},
\]
where $r(\alpha)=\min(\frac{1}{4},\frac{\alpha-1}{2\alpha
})$.
\end{lemma}

The next lemma extends the result on the growth rate of Wiener
processes in Lemma~\ref{lemma32} to the sequence of approximating
processes~(\ref{32}).

\begin{lemma}\label{lemma38} Suppose that \textup{(A2)} holds and that $p\in
(0,1]$. Then, as $k\rightarrow\infty$,
\[
\|W_{m,n}^0\|_{\infty} = \mathcal{O}\bigl(\sqrt{\log k}\bigr)
\qquad\mbox{a.s.}
\]
\end{lemma}

Having established the existence of a sequence $W^0_{m,n}$ of Gaussian
processes approximating $W_{m,n}$ and studied the behavior of $\mathcal
{F}^{-1}_{m,n}$, we may provide a sequence of Gaussian processes
approximating $U_{m,n}$. Define $U^0_{m,n}=\mathcal{F}^{-1}(W^0_{m,n})$
for each $m$ and $n$. Since $\mathcal{F}^{-1}$ is a bounded linear
operator, $U^0_{m,n}$ forms a sequence of Gaussian processes.
%
\begin{theorem}\label{theorem31}
Suppose that \textup{(A1)--(A4)} hold. Then, on the probability space on which
$W_{m,n}$ and $W^0_{m,n}$ are defined, we have that, as $k \rightarrow
\infty$,
\[
\| U_{m,n} - U_{m,n}^0 \|_{[0,\tau-\varepsilon_{m,n}]}
= \mathcal{O}\bigl(\varepsilon_{m,n}(\log k)^{{3}/{2}}\sqrt{\log
\log
k}\bigr) \qquad\mbox{a.s.},
\]
where $\varepsilon_{m,n}=k^{-{r(\alpha)}/({\beta+1})}$ and $r(\alpha
)=\min(\frac{1}{4},\frac{\alpha-1}{2\alpha})$.
\end{theorem}

Theorem~\ref{theorem31} will be crucial in our study of the
asymptotic properties of kernel density estimators of $g$, the density
associated to $G$, in Sections~\ref{kernel} and~\ref{ISE}. Other
applications of Theorem~\ref{theorem31} include oscillation moduli
and laws of the iterated logarithm; see~\cite{csorgo1984kmt}.

\subsection{Global modulus of continuity}

In order to describe the asymptotic properties of the kernel density
estimators of $g$ via the above strong approximation, we must establish
the global modulus of continuity of the approximating
process~$U_{m,n}^{0}$.

In the sequel, we say that the distribution $G$ satisfies assumption
(A5) if its density $g$ is differentiable,\vadjust{\goodbreak} and that a sequence
of bandwidths $h_{m,n}$ satisfies assumption (B1) if:
\begin{longlist}[(2)]
\item[(1)] $mh_{m,n} \to\infty$ and $\log h_{m,n}/\log\log m \to
-\infty$ as $m,n \to\infty$;
\item[(2)] $\sqrt{\log n}h_{m,n} \to0$ and $\sqrt{\log m}h_{m,n}\to0$
as $m,n \to\infty$.
\end{longlist}
%
\begin{theorem}\label{theoremG1}
Suppose that \textup{(A1)--(A5)} hold, and that the sequence $h_{m,n}$
satisfies \textup{(B1)}. Then, for any $\eta$ in $(0,\tau)$, we have
that, as $k \to \infty$,
\[
\sup_{0\leq t \leq\tau- \eta}\sup_{0\leq s \leq h_{m,n}
}|U_{m,n}^0(t+s) - U_{m,n}^0(t)| = \mathcal{O}\bigl(\sqrt{h_{m,n}
\log
(1/h_{m,n})}\bigr) \qquad\mbox{a.s.}\vspace*{-3pt}
\]
\end{theorem}


\section{Asymptotic behavior of kernel density estimators}
\label{kernel}

Consider the kernel density estimator $\hat{g}_m$ of a univariate
density $g$ introduced by~\cite{rosenblatt1956ams},
%
\begin{equation}\label{kernelform}
\hat{g}_m(t)=\frac{1}{h_m}\int_{0}^{\infty}{K\biggl(\frac
{t-s}{h_m}\biggr)\,d\hat{G}_m(s)},
\end{equation}
where $X_1,\ldots,X_m$ are independent observations from $g$, $K$ is
some kernel function, $h_m$ some bandwidth and $\hat{G}_m$ the
empirical distribution function. The weak and strong uniform
consistency of $\hat{g}_m$ was addressed in
\cite{nadaraya1965theoryprob,schuster1969ams} and~\cite
{vanryzin1969ams}, among others. To ensure strong uniform consistency,
these authors required that $\sum_m \exp(-cm{h_m}^2) < \infty$ for
each $c>0$. Silverman~\cite{silverman1978annals} established the
strong uniform
consistency of $\hat{g}_m$ under weaker assumptions using the KMT
strong approximation technique. When the observations are subject to random
right-censoring, Blum and Susarla~\cite{blum1980multanal} proposed
estimating $g$ by the estimator in~(\ref{kernelform}), replacing $\hat
{G}_m$ by the Kaplan--Meier estimator of $G$. The properties of the
resulting estimator were examined in
\cite{blum1980multanal,foldes1981pmh} and~\cite{mielniczuk1986annals},
among others.

To estimate the density function $g$ under multiplicative censoring, we
consider a sequence of kernel density estimators $\hat{g}_{m,n}$,
defined as
%
\begin{equation}\label{40}
\hat{g}_{m,n}(t) = \frac{1}{h_{m,n}}\int_0^\infty
K\biggl(\frac{t-s}{h_{m,n}}\biggr)\,d\hat{G}(s),
\end{equation}
where $\hat{G}$ is, as before, a solution of the nonparametric score
equation based on $(\mathcal{X}_m,\mathcal{Y}_n)$.

We introduce an additional set of assumptions to be used in the sequel.
The sequence of bandwidths $h_{m,n}$ is said to satisfy assumption
(B2) if
\[
\lim_{k\rightarrow\infty}\frac{\varepsilon_{m,n}(\log k)^{{3}/{2}}
\sqrt{\log\log k}}{\sqrt{k}h_{m,n}}=0.
\]
We say that a kernel function $K$ satisfies assumption (K1) if:
\begin{longlist}[(2)]
\item[(1)] $K$ has total variation $V_K<\infty$;
\item[(2)] $K$ is supported on $(-1,1)$;
\item[(3)] $K$ is continuous;
\item[(4)] $\int{K(u)\,du}=1$.
\end{longlist}
Further, we say that it satisfies assumption (K2) if $\int{uK(u)\,du}=0$.\vadjust{\goodbreak}

\subsection{Strong uniform consistency}\label{subsection41}

Denote by $g_{m,n}$ the kernel smoothing of $g$ based on $G$; that is,
write
\[
g_{m,n}(t)=\frac{1}{h_{m,n}}\int_{0}^{\infty}{K\biggl(\frac
{t-s}{h_{m,n}}\biggr)\,dG(s)}.
\]
%
\begin{lemma}\label{lemma412} Suppose that \textup{(A1)--(A5)} hold, and that
$h_{m,n}$ is a sequence of positive bandwidths tending to 0 as
$k\rightarrow\infty$ and satisfying \textup{(B1)} and \textup{(B2)}.
Suppose also that the kernel function $K$ satisfies \textup{(K1)}. Then, for any
$\eta$ in $(0,\tau )$, we have that
\[
{\lim_{k \rightarrow\infty}} \|\hat{g}_{m,n}-g_{m,n}\|
_{[0,\tau-\eta]}= 0 \qquad\mbox{a.s.}
\]
\end{lemma}
%
\begin{theorem}\label{theorem411} Suppose that \textup{(A1)--(A5)} hold, and
that $h_{m,n}$ is a sequence of positive bandwidths tending to 0 as
$k\rightarrow\infty$ and satisfying \textup{(B1)} and \textup{(B2)}.
Suppose also that the kernel function $K$ satisfies \textup{(K1)}. Then, for any
$\eta$ in $(0,\tau )$, we have that
\[
\lim_{k \rightarrow\infty} \|\hat{g}_{m,n}- g\|
_{[0,\tau
-\eta]}= 0 \qquad\mbox{a.s.}
\]
\end{theorem}



\subsection{Strong uniform approximation of the empirical density
process} By Theorems~\ref{theorem31}
and~\ref{theorem411}, we can find a sequence of Gaussian
processes that strongly and uniformly approximates the empirical density
process. Let $K$ be an arbitrary density function, and define
\[
\varphi_{m,n}(t,s) = \frac{1}{h_{m,n}} K\biggl(\frac
{t-s}{h_{m,n}}
\biggr).
\]
Denoting by $\mathrm{v}_{s} [\varphi_{m,n}(t, s)]$ the total
variation of $\varphi_{m,n}(t, \cdot)$ for fixed $t$, we refer to the
uniform total variation $\sup_t \mathrm{v}_{s} [\varphi_{m,n}(t,
s)]$ by $V_{m,n}$.
%
\begin{theorem}\label{theorem421}
Suppose that \textup{(A1)--(A5)} hold, and that $h_{m,n}$ is a sequence
of positive bandwidths tending to 0 as $k\rightarrow\infty$ and
satisfying \textup{(B1)} and \textup{(B2)}. Suppose also that the
kernel function $K$ satisfies \textup{(K1)} and \textup{(K2)}, and that~$g$ has a bounded
second derivative. Then, for any $\eta$ in $(0,\tau)$, we have that
\begin{eqnarray*}
&&\bigl\|\sqrt{k}(\hat{g}_{m,n} - g)-\Gamma_{m,n}\bigr\|_{[0,\tau
-\eta]}\\
&&\qquad= \mathcal{O} \biggl( \frac{\varepsilon_{m,n}(\log k)^{
{3}/{2}}\sqrt
{\log\log k}}{h_{m,n}} + \sqrt{k} h_{m,n}^{2} \biggr) \qquad\mbox{a.s.},
\end{eqnarray*}
where we have defined $\Gamma_{m,n}(t)=\int_0^\infty
U_{m,n}^0(s)\varphi
_{m,n}(t,ds)$.
\end{theorem}
%
\begin{remark}\label{remarkthm4}
Theorem~\ref{theorem421} suggests that the optimal rate for the
above approximation is obtained by choosing $h_{m, n} \sim(\varepsilon
_{m,n}\sqrt{\log\log k/k})^{{1}/{3}}\sqrt{\log k}$.
\end{remark}

Theorem~\ref{theorem421} implies distributional results. The
linearization $\psi_U(s-uh,t-vh)-\psi_U(s,t)\sim h$ is useful here.
This result is not difficult\vadjust{\goodbreak} to show for $p>1/2$ using representations
of $\psi_U$ provided on page 1033 of~\cite{vardi1992lannals},
linearization techniques and the modulus of continuity of process $U$.
The case $p\leq1/2$ (i.e., heavy censoring) is more challenging, but
can be dealt with using~(\ref{constraint1}),~(\ref{constraint2}) and an
argument similar to that found in the proof of
Theorem~\ref{theoremG1}. Using Theorem~\ref{theorem421} and the above
linearization, we
may show that \mbox{$\sqrt{kh_{m,n}}(\hat{g}_{m,n}-g)$} is
asymptotically Gaussian with mean zero and covariance function $\sigma
_{g}$ estimated consistently by
\[
\hat{\sigma}_{g}(s,t)=h_{m,n}^{-1}\iint{\hat{\psi
}_U(s-uh_{m,n},t-vh_{m,n})\,dK(u)\,dK(v)}.
\]
%

\section{Integrated squared error of kernel density
estimators}\label{ISE}
A common measure of the global performance of an estimator $\hat
{g}_{m}$ of a density $g$ is its integrated square error (ISE), defined as
\[
\mathcal{E}_{m} = \int_{-\infty}^\infty{[\hat
{g}_m(s)-g(s)
]^2\,ds}.
\]
Use of the ISE is particularly pervasive in simulation studies aiming
to compare the performance of various density estimators. 
Minimization of the mean integrated square error (MISE) $
\mathrm{E}[\mathcal{E}_m]=\int_{-\infty}^\infty{\mathrm{E}
[g_m(s)-g(s)]^2\,ds}$ is often a guiding principle in the
construction of kernel density estimators. Steele~\cite{steele1978cjs}
identified the need to determine the relationship between various
measures of accuracy in density estimation. One such measure, the order
of $\mathcal{E}_{m}-\mathrm{E}(\mathcal{E}_{m})$, is particularly
important in statistics. Hall~\cite{hall1982stochproc} first began
addressing the issues raised in~\cite{steele1978cjs} by computing the
exact order of convergence of $\mathcal{E}_{m}-\mathrm{E}(\mathcal
{E}_{m})$ to zero using the strong approximation technique developed by
Koml\'{o}s, Major and Tusn\'{a}dy~\cite{komlos1975zwvg} for the
standard empirical process. Zhang~\cite{zhang1998stochproc} studied the
case of random right-censoring using the strong approximation technique
of~\cite{burke1981probtheory} and~\cite{burke1988probtheory}. In this
section, we consider the ISE $
\mathcal{E}_{m,n}$ of the kernel estimator $\hat{g}_{m,n}$ under
multiplicative censoring and derive its asymptotic expansion. 

\subsection{Asymptotic expansion of the integrated squared error}

In the remainder of the paper, we make use of the following
assumptions. We say that the kernel function $K$ satisfies assumption
(K3) if it has finite second moment $\sigma^2>0$ and is
differentiable. Further, we say that the density $g$ satisfies
assumption (A6) if it is twice continuously differentiable. Of
course, assumption (A6) implies assumption~(A5).
Finally, we say that the sequence of bandwidths $h_{m,n}$ satisfies
assumption (B3) if
\[
\lim_{k\rightarrow\infty}
\frac{\sqrt{\log k} (\log\log k)^{{1}/{6}}}{h_{m,n} k^
{{1}/({\delta(\beta)})}}=0,
\]
where $\delta(\beta) = 4 + 4\beta/(2\beta+ 3)$. In the sequel, we
write $\nu$ for $\sqrt{\int_{-1}^{1}{K^{2}(u)\,du}}$.

The ISE of $\hat{g}_{m,n}$ on the interval $[u_1,u_2]$ is defined as
\[
\mathcal{E}_{m,n}(u_1,u_2) = \int_{u_1}^{u_2}{[\hat
{g}_{m,n}(s) -
g(s)]^{2} \,ds}.
\]
Theorem~\ref{theoremI1} presents an asymptotic
expansion for $\mathcal{E}_{m,n}(0,\tau- \eta)$ for any $\eta$ in
$(0,\tau)$.

%
%
%
%
%
%
\begin{theorem}\label{theoremI1}
Suppose that \textup{(A1)--(A4)} and \textup{(A6)} hold with
$\alpha_0>4$ in~\textup{(A3)} and that $\alpha$ is chosen in $[2,\alpha_0/2)$.
Suppose that $h_{m,n}$ is a sequence of positive bandwidths satisfying
\textup{(B1)} and \textup{(B3)}, and that the kernel function~K satisfies \textup{(K1)--(K3)}.
Then, for any $\eta$ in $(0,\tau)$, we have that
\[
\mathcal{E}_{m,n}(0,\tau-\eta)= \frac{h_{m,n}^4\sigma^4}{4}\int
_{0}^{\tau-\eta}{[g''(s)]^2\,ds}+\frac{\nu^2}{h_{m,n}
kp} +
o_{p} \biggl(\frac{1}{kh_{m,n}}+h_{m,n}^{4} \biggr) .
\]
\end{theorem}

Theorem~\ref{theoremI1} suggests that $h_{m,n}$ should shrink at the
rate $k^{-{1}/({\zeta(\beta)})}$ modulo logarithmic terms, where
$\zeta(\beta) = \max(5,\delta(\beta))$. We note that
$\delta
(\beta) < 5$ when $\beta< 3/2$. Then, writing
$\|g''\|^2_{2, [0, \tau-\eta]} = \int_{0}^{\tau- \eta}{
[g''(s)]^2\,ds}$, Theorem~\ref{theoremI1} suggests that the bandwidth
\[
h^{\star}_{m,n} = \biggl( \frac{\nu^2}{k p \sigma^4 \|g''\|_{2, [0,
\tau
-\eta]}^{2}} \biggr)^{1/5}
\]
minimizes the order of the integrated squared error, a direct
generalization of the reference rule for uncensored data alone, which
we recover for $p=1$ and $k=m$. Of course, in practice, this bandwidth
is unknown; instead, we may substitute $g''$ by some estimate $\hat
{g}''$, and $p$ by $\hat{p}=m/k$. For example, a~reference rule based
on a Gamma approximation to $G$ is given by
%
\begin{equation}\label{prule}
\hat{h}^{\star}_{m,n}=2\hat{\beta} \biggl(\frac{\nu^2}{ m \sigma
^4}
\biggr)^{1/5},
\end{equation}
where\vspace*{1pt} $\hat{\beta}=\sum_{i=1}^{m}{X_i}/(4m)$ is the MLE of $\beta$
based on $\mathcal{X}_m$ and the model $G = G_\beta$, with $g_\beta(x)
= x^3 \exp(-x/\beta)/(6\beta^4)$ the density associated to $G_\beta$.
This distribution satisfies (A3) with $\alpha_0=4$ but is a limiting
case with respect to the stronger assumption made in Theorem~\ref
{theoremI1}. It was selected because it has the least smooth density
in the family of densities $\{g_{\alpha,\beta}(x) = x^{\alpha-1}
\exp
(-x/\beta)/[\Gamma(\alpha)\beta^{\alpha}]\dvtx\alpha\geq4\}$ with respect
to the $L_2$-norm of the second derivative of $g_{\alpha,\beta}$.
Alternatively, we may consider kernel smoothing of the uncensored
observations alone to obtain a nonparametric pilot estimate~$\hat{g}''$
of~$g''$. More robust but computationally intensive cross-validation
approaches, as in~\cite{marron1987annals}, may also be used for
bandwidth selection.

\begin{figure}

\includegraphics{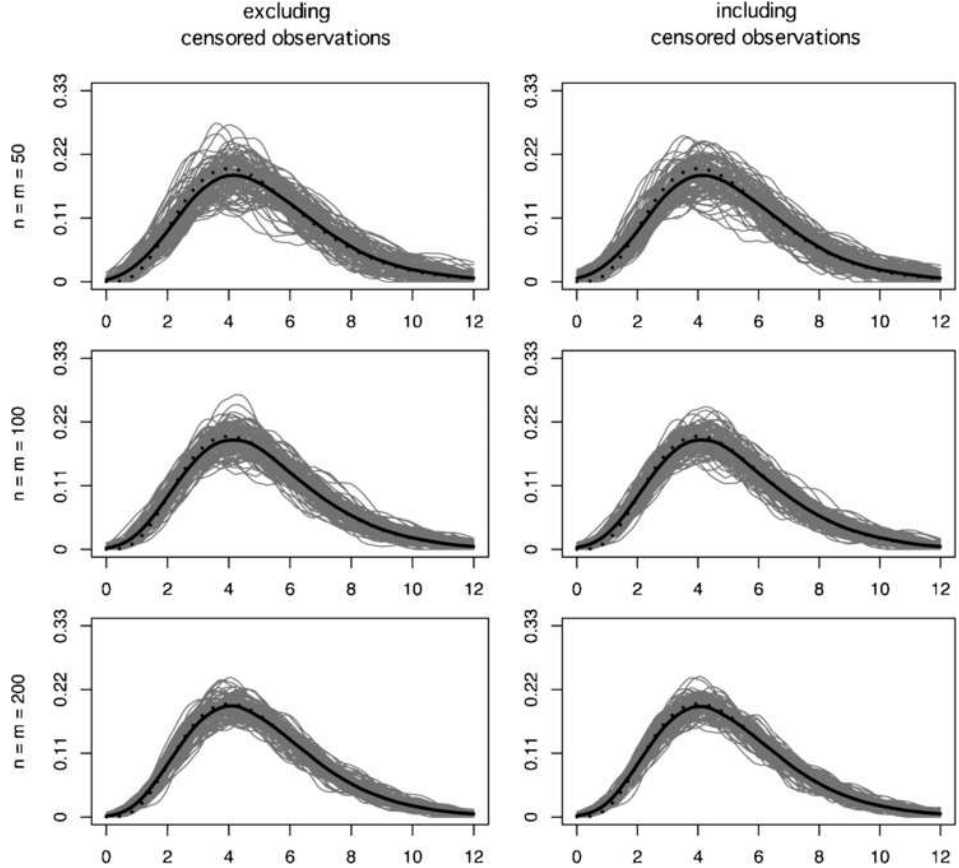}

\caption{Overlayed sample paths.} \label{plots}
\end{figure}

%

\subsection{Small-sample simulation results: Implementation and efficiency}

To provide some illustration of the behavior of the methods proposed,
we present below results from a preliminary small-sample\vadjust{\goodbreak} simulation
study. The objective was to graphically evaluate the general adequacy
of the estimators as well as to elucidate the potential contribution of
censored observations to overall estimation efficiency, both in small
samples. For this purpose, we considered data emanating from the
multiplicative censoring model, with underlying Gamma density function
$g_{\alpha}(x)=\Gamma(\alpha)^{-1}x^{\alpha-1}\exp(-x)\mathbb
{I}_{(0,\infty)}(x)$, various sample sizes and differing values of
parameter $\alpha$. We found the kernel density estimators proposed to
perform generally well. Figure~\ref{plots} presents 100 sample paths,
shown in grey, for various sample sizes and parameter value $\alpha=5$.
Plots in the first column were obtained by discarding all censored
observations and performing kernel density estimation using the
uncensored observations alone; all
observations were used in generating plots in the second column. The
pointwise average of the sample plots is shown in solid black, while
the true density is the dotted black curve depicted. The first, second
and third rows were generated from datasets of 100, 200 and 400 total
observations, respectively, with censored and uncensored observations
equally represented. In all cases, bandwidth values were automatically
selected using the $\Gamma(4,\beta)$ parametric reference rule~(\ref
{prule}). The Epanechnikov kernel $K(x)=\frac{3}{4}(1-x^2)\mathbb
{I}_{(-1,1)}(x)$ was used throughout. From these plots, we notice that
use of the full sample leads to a decrease in variability throughout
the support. Our empirical findings suggest that this cumulates to a
substantial decrease in integrated squared error. Table~\ref{table1}
reports estimates and associated 95\% confidence intervals for the mean
relative difference in ISE, defined as $(\mbox{ISE}_0-\mbox
{ISE}_1)/\mbox{ISE}_1$, obtained from a simulation of 500 datasets,
where $\mbox{ISE}_0$ and $\mbox{ISE}_1$ are the integrated squared
errors associated with the use of the uncensored subsample and of the
full sample, respectively. These values describe the mean percent
increase in $\mbox{ISE}$ from discarding the censored subsample, for
various sample sizes and parameter values.

\begin{table}
\caption{Average percent increase in ISE and 95\% CIs using $\Gamma
(4,\beta)$ parametric reference rule}\label{table1}
\begin{tabular*}{\tablewidth}{@{\extracolsep{\fill}}lcccc@{}}
\hline
\textbf{Sample size} & \multicolumn{1}{c}{$\bolds{\alpha=3}$} &
\multicolumn{1}{c}{$\bolds{\alpha=4}$} & \multicolumn{1}{c}{$\bolds{\alpha=5}$}
& \multicolumn{1}{c@{}}{$\bolds{\alpha=6}$} \\
\hline
$50+50$ & $_{15.2}17.1{}_{19.0}$ & $_{11.3}13.4{}_{15.6}$ &
$_{16.4}18.4{}_{20.4}$ & $_{13.9}16.3{}_{18.7}$ \\
$100+100$ & $_{16.8}18.5{}_{20.2}$ & $_{14.1}15.8{}_{17.5}$ &
$_{13.6}15.3{}_{17.0}$ & \hphantom{$_{9}$}$_{9.9}11.6{}_{13.2}$ \\
$200+200$ & $_{13.2}14.7{}_{16.1}$ & $_{11.6}13.1{}_{14.7}$ &
$_{14.4}17.8{}_{21.2}$ & $_{18.4}22.6{}_{26.7}$ \\
\hline
\end{tabular*}
\end{table}

The relative performance of the estimators was found to be rather
insensitive to the proximity of the underlying distribution to the
parametric model specified in the reference rule used, with an average
increase in ISE of around 10--25\%, subsequent to discarding censored
observations, regardless of sample size and parameter value. Since the
performance of kernel density estimators hinges upon the performance of
the underlying estimator of the distribution function as well as the
adequacy of the bandwidth selection rule, gauging the contribution of
censored observations to overall estimation efficiency is complicated
by the layer of uncertainty associated to bandwidth selection. As such,
we have also conducted a simulation study, whereby, for each simulated
dataset, the bandwidth selected was that minimizing the observed ISE;
we refer to this rule as the optimal bandwidth selection rule. Of
course, such a rule can only be adopted in simulation settings, where
the true density function is known, and the ISE can be computed
directly. This approach provides, nonetheless, a clearer view of the
gains resulting from the inclusion of censored observations in the
estimation procedure. Table~\ref{table2} reports estimates of the mean
relative increase in ISE resulting from discarding all censored
observations along with associated 95\% confidence intervals. These
results seem to suggest that for small and moderate sample sizes, when
equal numbers of censored and uncensored observations are available,
ignoring censored observations leads to an increase in ISE of roughly
10--35\%, results consistent with those reported in Table~\ref{table1}.

\begin{table}
\caption{Average percent increase in ISE and 95\% CIs using optimal
bandwidth selection rule}\label{table2}
\begin{tabular*}{\tablewidth}{@{\extracolsep{\fill}}lcccc@{}}
\hline
\textbf{Sample size} & \multicolumn{1}{c}{$\bolds{\alpha=3}$} &
\multicolumn
{1}{c}{$\bolds{\alpha=4}$} & \multicolumn{1}{c}{$\bolds{\alpha=5}$}
& \multicolumn{1}{c@{}}{$\bolds{\alpha=6}$} \\
\hline
$50+50$ & \hphantom{$_{9}$}$_{ 9.6} 14.3{}_{19.0}$ & $_{10.9} 15.7{}_{20.5}$
& \hphantom{$_{9}$}$_{ 9.8} 14.8{}_{19.9}$ & $_{17.3} 32.5{}_{47.7}$\\
$100+100$ & $_{12.7} 16.3{}_{20.0}$ & $_{12.9} 17.3{}_{21.6}$ &
$_{11.1} 15.8{}_{20.5}$ & $_{12.4} 26.9{}_{41.3}$ \\
$200+200$ & $_{10.0} 13.8{}_{17.6}$ & $_{12.1} 17.2{}_{22.4}$ &
$_{14.3} 21.0{}_{27.7}$ & $_{16.9} 34.6{}_{52.3}$\\
\hline
\end{tabular*}
\end{table}

The above provides a glimpse of the contribution of the censored
observations in small and moderate samples. It suggests that these
observations provide nonnegligible information regarding the estimand
of interest. We may, however, also resort to asymptotic arguments to
motivate use of the full sample for the sake of efficiency. For any
given distribution function $H$, denote the integrated squared error by
\[
\mbox{ISE}(H,h;g)=\int\biggl[\frac{1}{h}\int K\biggl(\frac
{x-y}{h}
\biggr)\,dH(y)-g(y)\biggr]^2\,dy
\]
and define the optimal bandwidth $\lambda(H;g)$ as the minimizer of the
ISE with respect to the true density $g$, that is, $\lambda
(H;g)=\argmin
_{h>0} \mbox{ISE}(H,h;g)$. Let~$G_{m,n}$ be any consistent estimator of
$G$ based on $(\mathcal{X}_m,\mathcal{Y}_n)$. The optimal kernel
density estimator of $g$ based on $G_{m,n}$ is then $g^{\star
}_{m,n}=\omega(G_{m,n})$, where~$\omega$ is the operator defined
pointwise as
\[
\omega(H)(x)=\frac{1}{\lambda(H;g)}\int K
\biggl(\frac{x-u}{\lambda(H;g)}\biggr)\,dH(u).
\]
Since any solution $\hat{G}$ of the nonparametric score equation is
asymptotically efficient for $G$ (see~\cite{vardi1992lannals}), it is
possible to show, along the lines of Theorem 25.47
of~\cite{vandervaart2000}, that $\hat{g}_{m,n}^{\star}=\omega(\hat{G})$
is asymptotically efficient for $g=\omega(G)$. In particular, the
kernel density estimator using the empirical distribution function
based on uncensored observations alone cannot be expected to be
asymptotically efficient, given that the latter is itself not efficient
for $G$. It is thus clear that, barring additional complications linked
to bandwidth selection, use of the full sample is preferable to that of
the uncensored subsample alone.

\section{Length-biased sampling with right-censoring} \label{lbs-rc}

As discussed in the\break \hyperref[intro]{Introduction}, the likelihood of length-biased
right-censored data is a particular case of that exhibited in~(\ref
{mc-like}). The literature on length-biased sampling can be traced as
far back as~\cite{wicksell1925biometrika}, with important contributions
by Fisher~\cite{fisher1934ema}, Neyman~\cite{neyman1955science} and
Zelen~\cite{zelen1969biometrika} in medical applications, and by
Cox~\cite{cox1969ssp} in industrial applications.
The rigorous treatment of biased sampling was initiated in the 1980s by
Vardi~\cite{vardi1982annals,vardi1985annals}, and furthered by
Gill, Vardi and Wellner~\cite{gill1988lannals}, Vardi and\vadjust{\goodbreak}
Zhang~\cite{vardi1992lannals}, Bickel and Ritov~\cite{bickel1991lannals},
Gilbert~\cite{gilbert2000annals} and, more recently, by Asgharian, M'Lan and
Wolfson~\cite{asgharian2002jasa}, Asgharian and Wolfson~\cite
{asgharian2005annals} and Bergeron, Asgharian and Wolfson~\cite
{bergeron2008jasa}. The importance of biased sampling in medical
applications and prevalent cohort studies was re-emphasized by Cox and
Oakes~\cite{cox1984book}.

The lifetime data typically collected on a prevalent cohort consist of
triples $(A, R\wedge D, \Delta)$, where $A, R$
and $D$ are, respectively, the current-age, the residual lifetime and
the residual censoring time, while $\Delta= \mathbb{I}_{\{R \leq D\}}$
is the censoring indicator. Suppose that $D$ and $(A, R)$ are
independent. In one scenario considered in~\cite{asgharian2005annals},
all analyses are
carried out conditionally upon the proportion of uncensored
individuals, assumed fixed. As such, the observations are comprised of
\[
(A_{i}, R_{i}) \stackrel{\mathrm{i.i.d.}}{\sim} f_{A, R\mid\Delta= 1},\qquad
i=1,\ldots, m,
\]
and
\[
(A_{j}, D_{j}) \stackrel{\mathrm{i.i.d.}}{\sim} f_{A, D\mid
\Delta= 0},\qquad j=1, \ldots, n,
\]
where $f_{A,R}(a,r) = f_U(a+r)/\mu_U$
and $f_{U}$ is the probability density function associated to
%
\begin{equation}\label{u-lb}
F_U(t) = \int_{0}^{t}{s^{-1}\,dG(s)}
{\Big/} \int_{0}^{\infty}{s^{-1}\,dG(s)}.
\end{equation}

The conditional density functions above are explicitly given by
\[
f_{A,R\mid\Delta=1}(a,r) = \frac{1-F_{D}(r)}{p(a+r)}\,dG(a+r)
\]
and
\[
f_{A,D\mid\Delta=0}(a,d) = \frac{f_{D}(d)}{(1-p)}\int
_{a+d\leq z}
z^{-1}\,dG(z)
\]
for the uncensored and censored subjects, respectively. Here, $f_D$ and
$F_D$ are, respectively, the density and distribution functions
associated to the residual censoring random variable $D$, and
$p=\mathrm{pr}(\Delta=1)$ is the proportion of uncensored individuals. The full
likelihood
of $m$ uncensored and $n$ censored length-biased observations is thus
\begin{eqnarray*}
\mathcal{L} &=& \prod_{i=1}^{m}{\biggl[\frac
{1-F_{D}(r_{i})}{px_{i}}\,dG(x_{i})\biggr]}
\prod_{j=1}^{n}\biggl[\frac{f_{D}(d_{j})}{1-p}\int_{y_{j}\leq
z}z^{-1}\,dG(z)\biggr] \\
&\propto& \prod_{i=1}^{m}\,dG(x_{i}) \prod_{j=1}^{n}
\int_{y_{j}\leq z}z^{-1}\,dG(z).
\end{eqnarray*}

Denoting $G_{\ast}(t)=P(A+R\leq t \mid\Delta=1)$ and $F_{\ast
}(t)=P(A+D\leq t \mid\Delta= 0)$ with associated density functions
$g_{\ast}(t)$ and $f_{\ast}(t)$, we may verify that
\[
g_{\ast}(t)= \frac{g(t)}{pt} \int_{0}^{t}{[1-F_{D}(r)
]\,dr}\quad\mbox{and}\quad f_{\ast}(t) = \frac{f(t)F_{D}(t)}{1-p},\vadjust{\goodbreak}
\]
where $f(t)$ is given by~(\ref{resid-life}). Defining the operators
\begin{eqnarray*}
\mathcal{H}(u)(t) &=& \int_{0<x \leq t}\frac{g_{\ast}(x)}{g(x)}\,du(x),
\\
\mathcal{K}_{m,n}(u)(t) &=& \int_{0<y \leq t }y \biggl(\int_{y \leq
z}\frac
{u(z)}{z^{2}} \,dz \biggr)\,
d\biggl[\biggl(\frac{\hat{f}(t)}{\hat{f}(y)} - 1\biggr) \frac
{f_{\ast}(y)}
{f(y)} \biggr]
\end{eqnarray*}
and $\Psi_{m,n}=\hat{p}\mathcal{H}+(1-\hat{p})\mathcal{K}_{m,n}$,
Asgharian and Wolfson~\cite{asgharian2005annals} have derived, under
this scenario, the equation $\Psi_{m,n}(U_{m,n}) = W_{m,n}$,
where $W_{m,n}$ is obtained from~(\ref{e4}) by replacing $W_{X,m}$ and
$W_{Y,n}$ by the empirical processes $\sqrt{m} (G_{m} - G_{*})$ and
$\sqrt{n}(F_{n} - F_{*})$, respectively. Defining the limiting
operators\looseness=1
\[
\mathcal{K}(u)(t) = \int_{0<y \leq t }y \biggl(\int_{y \leq z}\frac
{u(z)}{z^{2}}\,dz\biggr)
\,d\biggl[\biggl(\frac{f(t)}{f(y)} - 1\biggr) \frac{f_{\ast}(y)}
{f(y)} \biggr]
\]\looseness=0
and $\Psi=p\mathcal{H}+(1-p)\mathcal{K}$, one can show that $\Psi
_{m,n}$ converges almost
surely to $\Psi$ in operator norm topology, and that $\Psi$ is bounded,
linear and has bounded inverse $\Psi^{-1}$ if $p>0.59$; see~\cite
{asgharian2005annals}.


As discussed in the \hyperref[intro]{Introduction}, when the
observation mechanism generates length-biased samples, it is often of
prime interest to make inference about $F_{U}$ and its density function
$f_{U}$. Substitution\vspace*{1pt} of $G$ by $\hat{G}$ in~(\ref{u-lb})
yields~$\hat{F}_{U}$, an asymptotically efficient estimator of $F_U$. The
asymptotic properties of $Z_{m,n}=\sqrt{k}(\hat{F}_U-F_U)$ may be
studied via its relation to $U_{m,n}$. Indeed, defining $L_{t}(x) =
x^{-1}[\mathbb {I}_{[0,t]}(x) - F_{U}(t)]$, we may write
\[
\hat{F}_{U}(s)-F_{U}(s) = \frac{\int_{0}^{\infty} L_{s}(x) \,d[\hat{G}(x)
- G(x)]}
{\int_{0}^{\infty} x^{-1}\,d \hat{G}(x)},
\]
from which we have that $Z_{m,n}=\int_{0}^{\infty}{L_s(x)\,dU_{m,n}(x)}
/\int_{0}^{\infty}{x^{-1}\,d\hat{G}(x)} $. Defining the operator
$\mathscr {L}(g)(t) = \mu_{U}^{-1} \int_{0}^{\infty} L_{t}(x) \,dg(x)$,
we note that if there exists some $\gamma_0 >0$ such that $G(\gamma_0)=
0$ (in which case $G$ is said to satisfy assumption $\gamma$), the
operator $\mathscr{L}$ is bounded. Consequently, Theorems
\ref{theorem31}--\ref{theoremI1} hold when making inference about
$F_{U}$ and its density function $f_{U}$.



Under the additional assumption that the residual censoring
distribution does not have a point-mass at zero, it is possible to
provide an explicit distributional result for the empirical density
process arising from kernel density\vspace*{1pt} estimation.
Specifically, we have that the empirical density process
$\sqrt{kh_{m,n}}(\hat{f}_U-f_U)$ is asymptotically Gaussian with mean
zero and covariance function $\sigma _{f_U}$ estimated consistently by
\[
\hat{\sigma}_{f_U}(s,t)=h_{m,n}^{-1}\iint{\hat{\psi
}_Z(s-uh_{m,n},t-vh_{m,n})\,dK(u)\,dK(v)},
\]
where $\hat{\psi}_Z$ is a consistent estimator of the asymptotic
covariance function $\psi_{Z}$ associated to the sequence of processes
$Z_{m,n}$. For example, we may take
\[
\hat{\psi}_Z(s,t)=\hat{\mu}_U^{-2}\iint{\hat{\psi}_U(x,y)\,d\hat
{L}_s(x)\,d\hat{L}_t(y)},
\]
where $\hat{\mu}_U=\int_{0}^{\infty}{z^{-1}\,d\hat{G}(z)}$, $\hat
{L}_u(z)=z^{-1}[\mathbb{I}_{[0,u]}(z)-\hat{F}_U(z)]$ and $\hat{\psi}_U$
is a consistent estimator of the covariance function
$\psi_U$ of process $U$.
Since for $s\leq t$ we may write $\psi_U(s,t)$ as
\begin{eqnarray*}
\hspace*{-4.5pt}&&
p\biggl\{\int_{0}^{s}{[\beta(x)]^2\,dG_{*}(x)}-\biggl[\int
_{0}^{s}{\beta
(x)\,dG_{*}(x)}\int_{0}^{t} {\beta(x)\,dG_{*}(x)}\biggr]\biggr\}\\
\hspace*{-4.5pt}&&\quad{}
+
(1-p)\int_{0}^{t}\int_{0}^{s}f(x)f(y)\biggl\{e(x\wedge y)\\
\hspace*{-4.5pt}&&\hspace*{136.5pt}{}+h(x\wedge
y)\biggl[\frac{1}{f(x\vee y)}-\frac{1}{f(x\wedge y)}\biggr]
\\
\hspace*{-4.5pt}&&\hspace*{241pt}{}-h(x)h(y)\biggr\}\,d\zeta(x)\,d\zeta(y),
\end{eqnarray*}
where we have defined $\zeta(x)=g(x)[pg_{*}(x)]^{-1}$,
$h(x)=\int_{0}^{x}{F_{*}(y)\,d[1/f(y)]}$ and $e(x)=2\int
_{0}^{x}{h(y)\,d[1/f(y)]}$, consistent estimation of $\psi_U$ is possible
by substitution of appropriate empirical counterparts into the above.

Assumption $\gamma$ imposed on $G$ may seem restrictive, but
nonetheless holds in many industrial and medical applications. The case
of survival with dementia, studied in~\cite{asgharian2002jasa}
and~\cite{wolfson2001nejm}, is an example of such. It is possible to relax
this requirement by imposing that $G$ and $F_D$ vanish at zero at a
super-polynomial rate, that is, by assuming that $G(t)$ and $F_D(t)$
are $o(t^r)$ as $t\rightarrow0$ for each $r>0$. While preserving all
results pertaining to $G$, this relaxation does not directly preserve
those pertaining to $F_U$. The unboundedness of $\mathscr{L}$ is
problematic, although an application of Tikhonov's regularization
method may help in circumventing this problem. This has been explored
by Carroll, Rooij and Ruymgaart~\cite{carroll1991aam}, although not
from the perspective of strong approximations.

%

\section{Closing remarks} \label{conclusion}

(1) For distributions with a lighter left tail ($\alpha_0 > 2$) and
heavier right tail (small $\beta$), the rate obtained for the strong
approximation of $U_{m,n}$ is close to $k^{-1/4}$ modulo logarithmic
terms. It is unclear whether it is possible to achieve better rates; if
so, different techniques would necessarily be needed to control the
rate of $\mathcal{I}_5$ in Lemma~\ref{new4}, as the best achievable
rate for $\mathcal{I}_5$ using approximations by Bernstein polynomials
is $k^{-1/4}$. As for assumption (B2) on the bandwidth required to
establish Theorem~\ref{theorem411}, the $k^{-1/4}$ rate in the
strong approximation roughly translates into the bandwidth condition
$(\log
k)^2/(k^{3/4}h_{m,n}) \to0$ when we further replace the iterated
logarithmic term by a logarithmic term. This is in contrast to $\log k
/(kh_{m,n})\rightarrow0$ obtained in~\cite{silverman1978annals}, in
the case of uncensored observations alone. Likewise, the rate given in
Remark~\ref{remarkthm4}, after Theorem~\ref{theorem421}, is roughly
$h_{m,n} \sim(\log k)^{2/3}/k^{1/4}$.

(2) The theory presented in this paper requires that $\hat
{p}\rightarrow p\in(0,1]$. The case $p=0$ may itself be of interest. On
one hand, if $\hat{p}=0$ for each $k$, then all observations are
multiplicatively censored; this has been studied by Groeneboom~\cite
{groeneboom1985berkconf}, among others. On the other hand, if $\hat
{p}>0$ for each $k$, the methods developed in this paper may be adapted
as long as $\hat{p}$ does not vanish too rapidly. Specifically, we may
redefine $\mathcal{F}_{m,n}=\hat{p}\mathcal{I}+(1-\hat{p})\mathcal{G}$
and
\[
W^0_{m,n}(s) = \sqrt{\hat{p}}B_{X,m}(G(s))+\sqrt{1-\hat
{p}}f(s)\int_{0< y \leq
s}B_{Y,n}(F(y))\,d\biggl[\frac{1}{f(y)}\biggr].
\]
Suppose that $\hat{p}^{-2}$ is $\mathcal{O}(\upsilon_k)$ for some
sequence of positive real numbers $\upsilon_k$ tending to infinity.
%
%
Then the strong approximation holds, with $U^0_{m,n}$ redefined as the
Gaussian process $\mathcal{F}_{m,n}^{-1}(W^0_{m,n})$ and the rates
being multiplied by $\mathcal{O}(\upsilon_{k}^{2})$. Further, the rate
of the global modulus of continuity of $U^0_{m,n}$ is multiplied by
$\mathcal{O}(\upsilon_k)$.
This allows one to study the case $p=0$. This extension provides
insight into the leap between the square-root asymptotics in the
canonical multiplicative censoring setting and the cube-root
asymptotics for the Grenander estimator when only censored observations
are available.
%
%

\begin{appendix}\label{app}
\section*{Appendix: Proofs of main results}

\begin{pf*}{Proof of Claim~\ref{cla1}}
If the condition $\sum_{m,n}G(\gamma _{m,n})<\infty$ is satisfied, it
is an immediate consequence of Theorem 1 of Section 10.1 of
\cite{shorack1986} that $\mathrm{pr}(\min
(X_1,\ldots,X_m)\leq\gamma_{m,n} \mbox{ i.o.})=0$. Hence, almost
surely, we may find $m_0$ and $n_0\in\mathbb{N}$ such that, for each
$m\geq m_0$ and $n\geq n_0$, all uncensored observations
$x_1,\ldots,x_m$ are no smaller than $\gamma_{m,n}$. We restrict our
attention here to such sufficiently large $m$ and $n$. Define
$\delta_i=\mathbb{I}_{\{ x_1,\ldots,x_m\}}(t_i)$ for $i=1,\ldots,k$,
and write $r_0=\min\{i\dvtx t_i\geq \gamma_{m,n}\}$. By construction,
we must have that $\delta _1=\cdots=\delta _{r_0-1}=0$. Define the set
\[
\mathscr{D}=\Biggl\{(a_{r_0},a_{r_0+1},\ldots,a_{k})\dvtx 0\leq
a_{r_0},a_{r_0+1},\ldots,a_{k}\leq1,\sum_{i=r_0}^{k} a_i=1, a_{k}\geq
\frac{1}{k}\Biggr\},
\]
a bounded, closed and convex subset of $\mathbb{R}^{k-r_0+1}$. For
$i=r_0,\ldots,k$, define
\begin{eqnarray*}
\phi_i(a_{r_0},\ldots,a_{k}) &=& \delta_i\biggl(\frac{\hat
{p}}{m}
\biggr)+\frac{a_i}{t_i}\biggl(\frac{1-\hat{p}}{n}\biggr)\sum
_{j=1}^{i}{\frac
{1-\delta_j}{\sum_{q=\max(j,r_0)}^{k}a_q/t_q}}\\
&=&\frac{1}{k}\Biggl(\delta_i+\frac{a_i}{t_i}\sum_{j=1}^{i}{\frac
{1-\delta_j}{\sum_{q=\max(j,r_0)}^{k}a_q/t_q}}\Biggr)
\end{eqnarray*}
and $\phi=(\phi_{r_0},\ldots,\phi_{k})$. We note that $\phi$ is continuous
on $\mathscr{D}$. We want to show that $\phi(\mathscr{D})\subseteq
\mathscr{D}$. The fact that the image of $\mathscr{D}$ under $\phi_i$
is contained in $[0,1]$ for $i=r_0,\ldots,k$ is clear. That it is
contained in $[1/k,1]$ for $i=k$ is obvious if $\delta_k=1$. We assume
instead that $\delta_k=0$. Then, defining $\lambda_j=\sum_{q=\max
(j,r_0)}^{k-1}a_q/t_q\geq0$ for $j=1,\ldots,k-1$ and $\lambda
_k=a_k/t_k\geq0$, we observe that
\[
\frac{a_k}{t_k}\sum_{j=1}^{k}{\frac{1-\delta_j}{\sum_{q=\max
(j,r_0)}^{k}a_q/t_q}}=\lambda_k\Biggl(\sum_{j=1}^{k-1}{\frac
{1-\delta
_j}{\lambda_j+\lambda_k}}+\frac{1}{\lambda_k}\Biggr)\geq1,
\]
from which it follows that the image of $\mathscr{D}$ under $\phi_k$ is
contained in $[1/k,1]$ if $\delta_k=0$ as well. Finally,\vspace*{1pt}
we require the equality
$\sum_{i=r_0}^{k}{\phi_i(a_{r_0},\ldots,a_k)}=1$ to hold for any
$(a_{r_0},\ldots,a_k)\in\mathscr{D}$. This can be verified using that
\[
\sum_{i=r_0}^{k}\sum_{j=1}^{i}b_{ij}=\sum_{j=1}^{r_0-1}\sum
_{i=r_0}^{k}{b_{ij}}+\sum_{j=r_0}^{k}\sum_{i=j}^{k}{b_{ij}}
\]
for any array $b_{ij}$, where under the first sum on the right-hand
side, it holds that $\max(j,r_0)=r_0$, while under the second sum,
$\max
(j,r_0)=j$. We may thus use the Brouwer fixed point theorem (see, e.g.,
Proposition 2.6 on page~52 and Problem~6.7e on page 254 of~\cite
{zeidler1985}) to obtain that there exists some
$a^*=(a_{r_0}^*,\ldots,a_{k}^*)\in\mathscr{D}$ such that $\phi(a^*)=a^*$.
The distribution function
\[
\hat{G}^*(t)=\sum_{i=r_0}^{k}a^*_i\mathbb{I}_{[0,t]}(t_i)
\]
is a solution to equation~(\ref{sq1}) with zero mass below $\gamma_{m,n}$.
\end{pf*}
\begin{pf*}{Proof of Theorem~\ref{theorem31}}
Using Lemma~\ref{lemma31} and the boundedness of $\mathcal
{F}_{m,n}^{-1}$, we have for each $t\in[0,\tau-\varepsilon]$ that
\[
U_{m,n}(t) = \mathcal{F}_{m,n,\varepsilon}^{-1}(W_{m,n})(t) + \mathcal
{O}\bigl( \varepsilon\sqrt{\log\log k} \bigr) \qquad\mbox{a.s.}
\]
Similarly, using the definition of $U_{m,n}^{0}$, $W_{m,n}^{0}$, Lemma
\ref{lemma38} and the boundedness of~$\mathcal{F}^{-1}$, we have for
each $t\in[0,\tau-\varepsilon]$ that
\[
U_{m,n}^{0}(t) = \mathcal{F}_{\varepsilon}^{-1}(W_{m,n}^{0})(t) +
\mathcal{O}\bigl( \varepsilon\sqrt{\log k} \bigr) \qquad\mbox{a.s.}
\]
The result follows from Lemmas~\ref{lemma39},~\ref{lemma35} and~\ref
{lemma38} and the inequality
\begin{eqnarray*}
\|U_{m,n}-U_{m,n}^0\|_{[0,\tau-\varepsilon]}
&=& \|{\mathcal{F}}_{m,n}^{-1}(W_{m,n})-{\mathcal
{F}}^{-1}(W_{m,n}^0)\|_{[0,\tau- \varepsilon]}
\\
&\leq&\| {\mathcal{F}}_{m,n, \varepsilon}^{-1} \|\|
W_{m,n} - W_{m,n}^0 \|_{[0,\tau-\varepsilon]}
\\
&&{}+\| \mathcal{F}_{m,n, \varepsilon}^{-1}
- \mathcal{F}^{-1}_{\varepsilon} \|\| W_{m,n}^{0} \|
_{[0,\tau-\varepsilon]} \\
&&{}+\mathcal{O}\bigl(\varepsilon\sqrt{\log k}
\bigr) \qquad\mbox{a.s.}
\end{eqnarray*}
We therefore find that
\begin{eqnarray*}
\|U_{m,n}-U_{m,n}^0\|_{[0,\tau-\varepsilon]}
& \leq&
\mathcal{O}\biggl(\frac{k^{-r(\alpha)}\sqrt{\log k}(\log\log
k)^{
{1}/{4}}}{f(\tau-\varepsilon)}\biggr) \\
&&{}
+\mathcal{O}\biggl(\frac{k^{-({\alpha-1})/({2\alpha})}\log k\sqrt
{\log
\log k}}{f(\tau-\varepsilon)} \biggr) \mathcal{O}\bigl(\sqrt{\log
k}
\bigr)\\
&&{}+\mathcal{O}\bigl(\varepsilon\sqrt{\log k}
\bigr) \qquad\mbox{a.s.}
\end{eqnarray*}
The use of Lemma~\ref{lemma39} was justified by the fact that
$W^0_{m,n}$ is almost surely continuous. Since (A4) implies that
$f(\tau
-u)\sim u^{\beta}$ for $u$ small, the above bound has least order,
modulo logarithmic terms, for $\varepsilon=\varepsilon_{m,n}$.
\end{pf*}
\begin{pf*}{Proof of Theorem~\ref{theoremG1}}
Let $t\in[0,\tau- \eta]$ and $s\in
[0,h_{m,n}]$. By definition~(\ref{32}), linearity of $\mathcal
{F}^{-1}$ and the triangle inequality, we have that
%
\begin{equation}\label{boundonU0}
|U^0_{m,n}(t+s)-U^0_{m,n}(t)
| \leq I_m(s,t) + J_n(s,t),
\end{equation}
where we define
\begin{eqnarray*}
I_m(s,t) &=& |\mathcal{F}^{-1}(B_{X,m}\circ G)(t+s)-{ \mathcal
{F}^{-1}(B_{X,m}\circ G)(t)}|,
\\
J_n(s,t) &=& | \mathcal{F}^{-1}(\mathcal{H}_n)(t+s) - \mathcal
{F}^{-1}(\mathcal{H}_n)(t)|
\end{eqnarray*}
and
\[
\mathcal{H}_n(t)=f(t)\int_{0< y\leq t}B_{Y,n}(F(y))\,d\biggl[\frac
{1}{f(y)}\biggr].
\]

We first study $I_m(s,t)$. Writing $\varsigma(u)(\cdot)=\int
_{0}^{\infty
}{K(\cdot,x)u(x)\,dx}$ and noting that $\int_{0}^{\infty}{\mathcal
{A}_0(\cdot,x)u(x)\,dx}=\mathcal{G}(u)(\cdot)$ for each $u$, equations
(\ref{constraint1}) and~(\ref{constraint2}) imply that $\varsigma
(u)\equiv-(1-p)\mathcal{G}(u+p\varsigma(u))/p^2$. It follows from (A5)
that $M_1=\|f'\|_{[0,\tau]}<\infty$. We find that
\begin{eqnarray*}
&&
|\mathcal{G}(w)(t+s)-\mathcal{G}(w)(s)|
\\
&&\qquad\leq|f(t+s)-f(t)|\biggl|\int_{0<y\leq t}y\biggl(\int_{y
\leq z} \frac{w(z)}{z^{2}} \,dz \biggr) \,d
\biggl[\frac{1}{f(y)} \biggr]\biggr|\\
&&\qquad\quad{}+|f(t)|\biggl|\int_{t<y\leq t+s}y\biggl(\int_{y
\leq z} \frac{w(z)}{z^{2}} \,dz \biggr) \,d
\biggl[\frac{1}{f(y)} \biggr]\biggr|\\
&&\qquad\leq|f(t+s)-f(t)|\biggl[\frac{1}{f(t)}-\frac{1}{f(0)}\biggr]\|w\|
_{[0,\tau]}\\
&&\qquad\quad{}+|f(t+s)|\biggl[\frac{1}{f(t+s)} - \frac
{1}{f(t)}
\biggr]\|w\|_{[0,\tau]}\\
&&\qquad=|f(t+s)-f(t)|\biggl[\frac{f(0)-f(t)}{f(0)f(t)}+\frac
{1}{f(t)}\biggr]\|
w\|_{[0,\tau]}\\
&&\qquad\leq \frac{2}{f(\tau-\eta)}\|w\|_{[0,\tau]}|f(t+s)-f(t)|\\
&&\qquad \leq
\frac
{2M_1s}{f(\tau-\eta)}\|w\|_{[0,\tau]}
\end{eqnarray*}
from which it follows, using~(\ref{represent}), that
%
\begin{eqnarray}\label{nu}
|\varsigma(u)(t+s)-\varsigma(u)(t)| &\leq& \frac
{2(1-p)M_1s}{p^2f(\tau-\eta)}\|u+p\varsigma(u)\|_{[0,\tau
]}\nonumber\\[-8pt]\\[-8pt]
&\leq& \frac{2(1-p)M_1s}{p^2f(\tau-\eta)}(2+p\|\mathcal
{F}^{-1}\|
)\|u\|_{[0, \tau]}.\nonumber
\end{eqnarray}
Using~(\ref{represent}) once more, we then have that
\begin{eqnarray*}
&&
\sup_{0\leq t\leq\tau-\eta}\sup_{0\leq s\leq
h_{m,n}}I_{m}(s,t)
\\
&&\qquad\leq p^{-1}\sup_{0\leq t\leq\tau-\eta}\sup
_{0\leq
s\leq h_{m,n}}\bigl|B_{X,m}\bigl(G(t+s)\bigr)-B_{X,m}(G(t))\bigr|\\
&&\qquad\quad{} +\sup_{0\leq t\leq\tau-\eta}\sup_{0\leq s\leq
h_{m,n}}|\varsigma(B_{X,m}\circ G)(t+s)-\varsigma(B_{X,m}\circ
G)(t)|.
\end{eqnarray*}
%
%
Using (A5), we may show, as in~\cite{mason1983probtheory} and~\cite
{shorack1986}, that
\[
{\sup_{0\leq x \leq a_\tau}\sup_{0\leq y \leq M_0 h_{m,n}}}
|\mathcal{W}_{
X,m}(x+y)-\mathcal{W}_{X,m}(x)| = \mathcal{O}\bigl(\sqrt{h_{m,n}
\log(1/h_{m,n} )}\bigr)
\]
almost surely, where $a_\tau= G(\tau- \eta)$, $M_0=\|g\|_{[0,\tau]}$,
and $\mathcal{W}_{X,m}$ is the Wiener process associated with
$B_{X,m}$; see Lemma 1.4.1 of~\cite{csorgo1981}. Hence, by an
application of the MVT, $B_{X,m}\circ G$ has modulus of continuity
\[
\mathcal{O}\bigl(\sqrt{h_{m,n} \log(1/h_{m,n} )}\bigr)
\]
as well. In view of~(\ref{nu}) and the fact that $\|B_{X,m}\circ G\|
_{[0,\tau]}$ is $\mathcal{O}(\sqrt{\log m})$ almost surely,
we have that
\[
\sup_{0\leq t\leq\tau-\eta}\sup_{0\leq s\leq h_{m,n}}
|\varsigma
(B_{X,m}\circ G)(t+s)-\varsigma(B_{X,m}\circ G)(t)| =
\mathcal
{O}\bigl(\sqrt{\log{m}}h_{m,n}\bigr)
\]
almost surely. It follows from the discussion above then that
%
\begin{equation}\label{G4}
\sup_{0\leq t \leq\tau- \eta}\sup_{0\leq s \leq h_{m,n} }
I_{m}(s,t)=\mathcal{O}\bigl(\sqrt{h_{m,n} \log(1/h_{m,n} )}
\bigr) \qquad\mbox{a.s.}
\end{equation}
We now turn to $J_n(s,t)$. Defining
\[
J'_n(s,t)=|f(t+s)-f(t)|\int_{0<y\leq t}{
|B_{Y,n}(F(y))|\,d\biggl[\frac{1}{f(y)}\biggr]}
\]
and
\[
J''_n(s,t)=|f(t+s)|\int_{t<y\leq
t+s}{
|B_{Y,n}(F(y))|\,d\biggl[\frac{1}{f(y)}\biggr]},
\]
we notice that $|\mathcal{H}_n(t+s)-\mathcal{H}_n(t)
|\leq
J'_n(s,t)+J''_n(s,t)$.
%
%
Using the MVT, we have that
\[
J'_n(s,t)\leq{\frac{M_1 s}{f(\tau-\eta)}\sup_{0\leq y\leq1}}{
|B_{Y,n}(y)|}
\]
and
\[
J''_n(s,t)\leq{\frac{f(t)-f(t+s)}{f(t)}\sup_{0\leq y\leq1}}{
|B_{Y,n}(y)|}\leq{\frac{M_1 s}{f(t)}\sup_{0\leq y\leq1}}{
|B_{Y,n}(y)|},
\]
so that $\sup_{0\leq t \leq\tau- \eta}\sup_{0\leq s \leq h_{m,n}
}J'_n(s,t)$, $\sup_{0\leq t \leq\tau- \eta}\sup_{0\leq s \leq h_{m,n}
}J''_n(s,t)$ and consequently ${\sup_{0\leq t \leq\tau- \eta}\sup
_{0\leq s \leq h_{m,n} }}|\mathcal{H}_n(t+s)-\mathcal{H}_n(t)|$ are
$\mathcal{O}(\sqrt{\log{n}}h_{m,n})$\vspace*{1pt} almost surely. Further,
using~(\ref{nu}), we have that
\begin{eqnarray*}
&&
{\sup_{0\leq t \leq\tau- \eta}\sup_{0\leq s \leq h_{m,n} }}
|\varsigma(\mathcal{H}_n)(t+s)-\varsigma(\mathcal{H}_n)(t)
|
\\
&&\qquad\leq\frac{2(1-p)M_1}{p^2f(\tau-\eta)} (2+p\cdot\|\mathcal
{F}^{-1}\|)\|\mathcal{H}_n\|_{[0, \tau]}h_{m,n} = \mathcal
{O}\bigl(\sqrt{\log n}h_{m,n}\bigr) \qquad\mbox{a.s.}
\end{eqnarray*}
so that\vspace*{1pt} $\sup_{0\leq t \leq\tau- \eta}\sup_{0\leq s \leq
h_{m,n}}J_n(s,t)=\mathcal{O}(\sqrt{\log n}h_{m,n})$ almost
surely using~(\ref{represent}). The theorem follows in view of this
last result,~(\ref{boundonU0}) and~(\ref{G4}).
\end{pf*}
\begin{pf*}{Proof of Theorem~\ref{theorem411}}
By the continuity (and hence uniform continuity) of $g$ on $[0,\tau]$,
the dominated convergence theorem may be used to show that
%
\begin{equation}\label{diffing2}
{\lim_{k \rightarrow\infty}} \sup_{0 \leq s \leq\tau- \eta}|g_{
m,n}(s) - g(s)|=0.
\end{equation}
The theorem follows immediately from Lemma~\ref{lemma412} and the
triangle inequality.
\end{pf*}
\begin{pf*}{Proof of Theorem~\ref{theorem421}}
By Theorem~\ref{theorem31} and integration
by parts, for any $t\in[0,\tau-\eta]$, we may write that
\begin{eqnarray*}
\hat{g}_{m,n}(t) - g(t)
&=& [\hat{g}_{m,n}(t) - g_{m,n}(t)] + [g_{m,n}(t) -
g(t)]\\
&=& \frac{1}{\sqrt{k}}\int_0^\infty{U_{m,n}^0(s)\,d\psi
_{m,n}(t,s)}\\
&&{}+\mathcal{O} \biggl(\frac{V_{m,n}\varepsilon_{m,n}(\log
k)^{{3}/{2}}\sqrt{\log\log k}}{\sqrt{k}} + \delta_{m,n}
\biggr) \qquad\mbox{a.s.},
\end{eqnarray*}
where $\delta_{m,n} = \sup_{0 \leq t \leq\tau-\eta}
|g_{m,n}(t)-g(t)|$. The result follows from
\cite{Rev72}, which shows that $\delta_{m,n}=\mathcal{O}
(h^2_{m,n}
)$ and $V_{m,n}=\mathcal{O}(1/h_{m,n})$.\vadjust{\goodbreak}
\end{pf*}
\begin{pf*}{Proof of Theorem~\ref{theoremI1}}
Since $g$ is twice continuously differentiable on $[0,\tau-\eta]$, we
may write that
$g_{m,n}(s) - g(s) = h_{m,n}^2 \sigma^2 g''(s)/2 + o
(h_{m,n}^2)$
uniformly in $s \in[0, \tau- \eta]$. Combining this expansion with
(S.1) in the proof of Lemma 7 (see supplementary material~\cite
{asgharian2012supp}) yields
\begin{eqnarray*}
\hat{g}_{m,n}(s)-g(s) &=& \biggl(\frac{h_{m,n}^2 \sigma^2}{2}
\biggr)g''(s) +\frac{\Upsilon_{m,n}(s,h_{m,n})}{\sqrt{k}h_{m,n}}
\\
&&{}
+\mathcal{O} \biggl(\frac{\varepsilon_{m,n}(\log
k)^{
{3}/{2}}\sqrt{\log\log k}}{\sqrt{k}h_{m,n}}\biggr)+o
(h_{m,n}^2
) \qquad\mbox{a.s.}
\end{eqnarray*}
uniformly in $s \in[0, \tau- \eta]$, where $\Upsilon
_{m,n}(s,h)=\int
_{-1}^{1}{U_{m,n}^0(s-uh)\,dK(u)}$. In view of~(\ref{represent}) and the
proof of Theorem~\ref{theoremG1}, we find that
\[
\Upsilon_{m,n}(s, h_{m,n}) = p^{-{1}/{2}} \int_{-1}^{1}
B_{X,m}\bigl(G(s-uh_{m,n})\bigr)\,dK(u) + \mathcal{O}\bigl( \sqrt{\log k}h_{m,n}
\bigr) \quad\mbox{a.s.}
\]
Further, using (B3) we may show, for $\alpha\geq2$, that
\[
\frac{\varepsilon_{m,n}(\log k)^{{3}/{2}}\sqrt{\log\log k}}{\sqrt
{k}h_{m,n}} = o (h_{m,n}^{2} )
\]
%
and therefore that
\[
\hat{g}_{m,n}(s)-g(s)=\biggl(\frac{h_{m,n}^2 \sigma^2}{2} \biggr)g''(s)
+ \frac{\int_{-1}^{1} B_{X,m}(G(s-uh))\,dK(u)}{\sqrt{pk}h_{m,n}} +
o(h_{m,n}^{2})
\]
almost surely. It then follows that $\mathcal{E}_{m,n}(0,\tau-\eta)$
may be written as
\begin{eqnarray*}
&&\frac{h_{m,n}^4 \sigma^4}{4} \int_{0}^{\tau- \eta} {\{
g''(s)
\}^2\,ds}
+\frac{_{\eta}P_{m,n}(h_{m,n})}{p k h_{m,n}^2} + \frac{\sigma^2
h_{m,n}{} _{\eta}Q_{m,n}(h_{m,n})}{\sqrt{pk}}
\\
&&\qquad{}+o(h_{m,n}^2) \biggl\{ o(h_{m,n}^2) + h_{m,n}^2 \sigma^2
\int
_{0}^{\tau- \eta}g''(s)\,ds + \frac{2 _{\eta}R_{m,n}(h_{m,n})}{\sqrt
{pk}h_{m,n}}\biggr\} \qquad\mbox{a.s.},
\end{eqnarray*}
where we have defined
\begin{eqnarray*}
_{\eta}P_{m,n}(h)&=&\int_{0}^{\tau-\eta} {\biggl[\int
_{-1}^{1}B_{X,m}\bigl(G(s-uh)\bigr)\,dK(u)\biggr]^2\,ds} ,\\
_{\eta}Q_{m,n}(h)&=&\int_{0}^{\tau-\eta}{g''(s) \biggl[\int
_{-1}^{1}B_{X,m}\bigl(G(s-uh)\bigr)\,dK(u)\biggr]\,ds}
\end{eqnarray*}
and
\[
_{\eta}R_{m,n}(h)=\int_{0}^{\tau-\eta} \biggl[\int
_{-1}^{1}B_{X,m}\bigl(G(s-uh)\bigr)\,dK(u)\biggr]\,ds.
\]
It follows\vspace*{1pt} from~\cite{hall1982stochproc} that $_{\eta}P_{m,n}(h) =
h_{m,n} \nu^{2} + o_{p}(h_{m,n})$, while $_{\eta}Q_{m,n}(h)$ and $
_{\eta}R_{m,n}(h)$ are both $o_{p}(\sqrt{h_{m,n}})$. We therefore
obtain that $\mathcal{E}_{m,n}(0,\tau-\eta)$ may be expressed as
\[
\frac{h_{m,n}^4 \sigma^4}{4} \int_{0}^{\tau- \eta} {
[g''(s)
]^2\,ds} + \frac{\nu^{2}}{pkh_{m,n}} + o_{p}\Biggl(\frac
{1}{kh_{m,n}}+h_{m,n}\sqrt{\frac{h_{m,n}}{k}}+h_{m,n}^{4}\Biggr).
\]
The result follows upon noticing that a term of order
$o_{p}(h_{m,n}^{3/2}/\sqrt{k})$ is dominated by any term of order
$o_{p}(h_{m,n}^{4})$.
\end{pf*}
\end{appendix}

\section*{Acknowledgments}

The authors thank the current and previous Editors, Professors Cai and
Silverman, the anonymous Associate Editor and the referees for their
exceptionally insightful and constructive
comments. 

\begin{supplement}[id=suppA]
\stitle{Additional technical details: Proof of lemmas\\}
\slink[doi]{10.1214/11-AOS954SUPP}
\sdatatype{.pdf}
\sfilename{aos954\_supp.pdf}
\sdescription{The proof of each lemma in the paper is provided in the
supplementary material.}
\end{supplement}


\printaddresses

\end{document}